\documentclass[12pt,draft]{amsart}

\makeatletter

\theoremstyle{plain}    
\newtheorem{thm}{Theorem}[section]
\numberwithin{equation}{section} %% Comment out for sequentially--numbered
\numberwithin{figure}{section} %% Comment out for sequentially--numbered
\theoremstyle{plain}    
\newtheorem*{thm*}{Theorem} 
\theoremstyle{plain}    
\theoremstyle{plain}    
\newtheorem{lem}[thm]{Lemma} %%Delete [thm] to re--start numbering
\theoremstyle{plain}    
\newtheorem{prop}[thm]{Proposition} %%Delete [thm] to re--start numbering
\theoremstyle{definition}
\newtheorem{defn}[thm]{Definition}
\theoremstyle{remark}
\newtheorem{rem}[thm]{Remark}
\theoremstyle{remark}

\theoremstyle{remark}    

\theoremstyle{remark}    

\theoremstyle{definition}  

\theoremstyle{remark}
  \newtheorem*{acknowledgement*}{Acknowledgement} 

\theoremstyle{plain}    

\theoremstyle{plain}    
\theoremstyle{plain}    
\theoremstyle{plain}    
\theoremstyle{definition}

\theoremstyle{remark}

\theoremstyle{remark}    

\theoremstyle{remark}    

\theoremstyle{plain}    

\usepackage{amssymb} \usepackage{amscd}

\makeatother

\begin{document}

\title[Quasidiagonality]{Quasidiagonality and the finite section method}

\author{Nathanial P. Brown}

\address{Department of Mathematics, Penn State University, State
College, PA 16802}

\email{nbrown@math.psu.edu}

\thanks{ Partially supported by an NSF Postdoctoral
Fellowship and DMS-0244807.}

\begin{abstract}
Quasidiagonal operators on a Hilbert space are a large and important
class (containing all self-adjoint operators for instance).  They are
also perfectly suited for study via the finite section method (a
particular Galerkin method).  Indeed, the very definition of
quasidiagonality yields finite sections with good convergence
properties.  Moreover, simple operator theory techniques yield
estimates on certain rates of convergence.  In the case of
quasidiagonal band operators both the finite sections and rates of
convergence are explicitly given.
\end{abstract}

\maketitle

%%%%%%%%%%%%%%%%%%%%%%%%%%%%%%%%%%%%%%%%%%%%%%%%%%%%%%%%%%%%%%%%%%%%%
\section{Introduction}

This paper follows up on an observation from \cite[Section
10]{brown:AFD}.  Namely we see what happens when one takes the general
theorems from \cite{HRS} and applies them to the particular class of
quasidiagonal operators.  The end results are quite satisfying both
from a theoretical point of view and a practical computational point
of view.  We will see that the finite section method works very well
on the class of quasidiagonal operators in general and some explicit
and reasonable rates of convergence can be worked out in some interesting 
cases.

Let $H$ be a separable, complex Hilbert space, $B(H)$ denote the
bounded, linear operators on $H$ and $T \in B(H)$ be given.  An
operator $T$ is called {\em quasidiagonal} if there exist finite rank
orthogonal projections $P_1, P_2, P_3 \ldots$ such that $ \| P_n (v) -
v \| \to 0$ for all vectors $v \in H$ and which asymptotically commute
with $T$; $$\| [T,P_n] \| = \| TP_n - P_n T \| \to 0.$$ This notion is
due to Halmos (cf.\ \cite{halmos}, \cite{brown:QDsurvey},
\cite{dvv:QDsurvey}) and is inspired by the classical Weyl-von Neumann
Theorem; if $S \in B(H)$ is a {\em self-adjoint} operator then $S = D
+ k$ where $D$ is a diagonal operator with respect to some orthonormal
basis of $H$ and $k$ is a compact operator.  Indeed, Halmos observed
that $T \in B(H)$ is quasidiagonal if and only if $T = B + k$ where
$B$ is a block diagonal operator with respect to some orthonormal
basis of $H$ and $k$ is compact. (This means there exist projections
as in the definition of quasidiagonality but with the stronger
properties that $P_n H \subset P_{n+1}H$ and $[P_n, B] = 0$ for all
$n$ -- writing $H = P_1 H \oplus (P_2 - P_1)H \oplus (P_3 - P_2)H
\oplus \cdots$, the matrix of $B$ with respect to this decomposition
is easily seen to be block diagonal.)  In particular, every
self-adjoint operator is quasidiagonal and, more generally, a result
of Berg implies that every normal operator is also quasidiagonal (cf.\
\cite[Corollary II.4.2]{davidson}).  The self-adjoints contain a large
set of interesting examples (Almost Mathieu operators, discretized
Hamiltonians, etc.) and there are plenty of interesting
non-self-adjoint, but still normal, examples like difference operators
coming from PDE's with constant coefficients.  There are even
interesting non-normal, but still quasidiagonal, operators such as
certain weighted shifts and Toeplitz operators (not all of these are
quasidiagonal but many examples are known to be quasidiagonal) to
which the ideas of this paper apply.

Section 2 of this note contains a rather extensive list of definitions
and recalls some of the theorems from \cite{HRS} that we will need.
There is nothing new in this section but we hope the unfamiliar reader
will find it a convenient reference.

Sections 3, 4, 5 and 6 are devoted to some general theory.  These
sections are the theoretical foundation on which the actual numerical
analysis will be based.  That is, we make some very general
observations which are of little computational value by themselves.
In some sense, one can regard these results as analogues of an
existence result for a large class of PDE's.  When faced with an
actual PDE it is, of course, good to know via some general theory that
a solution exists and even better if one knows it is unique.  However,
the real goal of numerical analysis is to find explicit algorithms
which are guaranteed to converge to a solution and, moreover, study
the rate of convergence.  With this analogy in mind one can regard
sections 3 through 6 as proving existence and uniqueness results as
well as giving general estimates on rates of convergence.  Then in
sections 7 and 8 some explicit algorithms are described which put the
general machine into motion in two special cases.  Below we give a
more detailed description of the contents of this paper.

In Section 3 we observe the {\em existence} of good finite sections
for quasidiagonal operators and the {\em uniqueness} of the limits of
things like norms, condition numbers, spectral quantities, etc.\ of
the corresponding finite sections.  

Section 4 gives some simple estimates on rates of convergence for the
finite section method.  Though the rates of convergence are not always
very good they are best possible (at this level of generality) as
shown by simple examples.  (However, in certain cases the rates of
convergence can be quite good.)  Section 5 simply takes the results of
sections 3 and 4 and specializes them to the cases of normal and
self-adjoint operators where the rates of convergence have cleaner
statements.  In section 6 we quickly point out that the resolution of
a well-known problem in operator theory (Herrero's approximation
problem) yields an improved  existence result for many quasidiagonal
operators.  Putting this result into practice in special cases seems
quite difficult but a related strategy, inspired by the observations
of section 6, will be the subject of another paper (cf.\
\cite{brown:PV}).

Sections 7 and 8 are devoted to explicit descriptions of how to
implement the quasidiagonal approach in the cases of unilateral and,
respectively, bilateral band operators.  In the unilateral case we use
canonical finite sections and translate our previous results into
computable error estimates.  In the bilateral case a famous technique
of Berg turns out to work beautifully and, again, computable rates of
convergence are given.

\begin{acknowledgement*}
I thank Bill Arveson and Marc Rieffel for sparking my interest in
numerical analysis and, more generally, for providing more support and
inspiration to me than they probably realize.  I am also indebted to
the authors of \cite{HRS} and \cite{bottcher} as they have made a
large body of work readily accessible to the C$^*$-community.
Finally, I thank my colleague Ludmil Zikatanov for numerous helpful
conversations.
\end{acknowledgement*}

%%%%%%%%%%%%%%%%%%%%%%%%%%%%%%%%%%%%%%%%%%%%%%%%%%%%%%%%%%%%%%%%%%%%%
\section{Notation, definitions and known results}

This section is devoted to setting notation, recalling a number of
definitions and stating some theorems from \cite{HRS} that will be
required.

Let $T \in B(H)$ be a bounded operator on some complex Hilbert space
$H$ and $P_1 \leq P_2 \leq \cdots$ be a {\em filtration} (i.e.\ each
$P_n = P_n^* = P_n^2$ is a finite rank, orthogonal projection, the
range of $P_n$ is contained in the range of $P_{n+1}$ and $\| P_n(v) -
v \| \to 0$ for all vectors $v \in H$).  If the rank of $P_n$ is some
integer $k(n)$ then we will often -- but not always -- identify $P_n
B(H) P_n$ with the complex $k(n) \times k(n)$ matrices $M_{k(n)}
({\mathbb C})$.  We will refer to the sequence of operators $(P_n
TP_n)$ as the {\em finite sections} (of $T$ with respect to $\{ P_n
\}$). The {\em finite section method} is simply to try and study the
(infinite dimensional) operator $T$ via its finite sections with
respect to a suitable filtration.  For example, can one recover the
spectrum of $T$ as some sort of limit of the spectra of (the finite
dimensional matrices) $P_n TP_n$?  Of course, the choice of filtration
is very important as it can easily happen that the finite sections do
not approximate the original operator in any reasonable sense.  A
classic example is that of the unilateral shift $S$ on $l^2({\mathbb
N})$.  With respect to the canonical basis, the matrix of $S$ is
below.
$$\begin{pmatrix}
0 & 0 & 0 & 0 & \cdots\\
1 & 0 & 0 & 0 & \cdots\\
0 & 1 & 0 & 0 & \cdots\\
0 & 0 & 1 & 0 & \cdots\\
\vdots & \vdots & \vdots & \vdots & \ddots
\end{pmatrix}$$
Cutting $S$ to the canonical filtration yields nilpotent matrices and
thus the spectra of the finite sections are all one point sets $\{
0\}$ which do not converge to the spectrum of $S$ (which is the unit
disc).

\begin{center}
{\em Some operator theory definitions.}
\end{center}

\begin{defn} Given an operator $T \in B(H)$:
\begin{enumerate}
\item For $\varepsilon > 0$, the $\varepsilon$-{\em pseudospectrum} of
$T$ is $$\sigma^{(\varepsilon)}(T) = \{ \lambda \in {\mathbb C}: \|
(\lambda - T)^{-1} \| \geq \frac{1}{\varepsilon} \},$$ where $\lambda
- T$ denotes the operator which maps $v \in B(H)$ to $\lambda v - Tv$
and we will adopt the convention that if an operator $X \in B(H)$ is
{\em not invertible} then $\| X^{-1}\| = \infty$.  In particular, the
usual spectrum $\sigma(T)$ is contained in $\sigma^{(\varepsilon)}(T)$
for every $\varepsilon > 0$ and, moreover, $$\sigma(T) =
\bigcap_{\varepsilon > 0} \sigma^{(\varepsilon)}(T).$$ See 
{\bf http://web.comlab.ox.ac.uk/projects/pseudospectra} for more,
including software for computing pseudospectra.

\item The {\em singular values} of $T$ are $\sigma_2(T) = \{ \lambda >
0: \lambda^2 \in \sigma(T^*T) \}.$

\item We let $C^*(T)$ denote the {\em unital, C$^*$-algebra generated
by $T$}.  It is the smallest unital C$^*$-subalgebra of $B(H)$ which
contains $T$ and can be realized as the norm closure of the set of all
polynomials (where $T^0$ is defined to be the identity operator) in
$T$ and $T^*$.  A {\em state} on $C^*(T)$ is a linear functional $\phi
: C^*(T) \to {\mathbb C}$ such that $\| \phi \| = \phi(I) = 1$, where
$I$ is the identity operator on $H$. The set of all states is denoted
$S(C^*(T))$.

\item The {\em spatial numerical range} of $T$ is $$SN(T) = \{
\langle Tv,v\rangle : v \in H, \| v \| = 1\}$$ and the {\em algebraic
numerical range} of $T$ is $$AN(T) = \{ \phi(T) : \phi \in
S(C^*(T))\}.$$ Note that $SN(T) \subset AN(T)$. 

\item $T$ is called {\em Moore-Penrose invertible} if there exists an
operator $T^+ \in B(H)$ such that $TT^+ T = T$, $T^+ T T^+ = T^+$,
$(TT^+)^* = TT^+$ and $(T^+ T)^* = T^+ T$.  If such $T^+$ exists then
it is unique and is called the Moore-Penrose inverse of $T$.  It turns
out that if $T^+$ exists then $T^+ \in C^*(T)$ (cf.\ \cite[Theorem
2.15]{HRS}) and hence, for an abstract C$^*$-algebra the notion of
Moore-Penrose invertible does not depend on the choice of faithful
representation.  Also, recall that $T^+$ determines {\em least square
solutions} to operator equations; if $y \in H$ then $x = T^+ y$ is the
least square solution of the equation $Tx = y$ (i.e.\ $x$ is the
unique vector such that $\| Tx - y \| = \inf \{ \| Tu - y\|: u \in
H\}$ and $\| x \| = \inf\{ \|u\|: \| Tu - y \| = \inf \{ \| Tu - y\|:
u \in H\} \}$).

\item If $T$ is invertible then the {\em condition number} of $T$ is
$$cond(T) = \| T\| \| T^{-1}\|,$$ while the {\em generalized condition
number} of a Moore-Penrose invertible operator $T$ is the number $\|
T\| \| T^{+}\|$.
\end{enumerate} 
\end{defn}

We now need to recall some definitions motivated by the finite section
method, but this will require introducing a bit more notation.  If
$P_1 \leq P_2 \leq \cdots$ is a filtration of $H$ and the rank of
$P_n$ is $k(n)$ then it is useful to build a C$^*$-universe where
finite sections can be treated as a single element.  Namely we let
$$\Pi M_{k(n)} ({\mathbb C}) = \{ (x_n): x_n \in M_{k(n)} ({\mathbb
C}), \sup \|x_n \| < \infty\}$$ and $$\oplus M_{k(n)} ({\mathbb C}) =
\{ (x_n) \in \Pi M_{k(n)} ({\mathbb C}): \lim \| x_n \| = 0 \}$$
denote the $l^{\infty}$ and, respectively, $c_0$ direct sum of the
$M_{k(n)} ({\mathbb C})$'s (where we identify $P_nB(H)P_n$ and
$M_{k(n)} ({\mathbb C})$).  Note that $\oplus M_{k(n)} ({\mathbb C})$
is a closed, two-sided ideal in $\Pi M_{k(n)} ({\mathbb C})$ and hence
there is a canonical quotient mapping $\pi : \Pi M_{k(n)} ({\mathbb
C}) \to \Pi M_{k(n)} ({\mathbb C}) / \oplus M_{k(n)} ({\mathbb C})$.
This quotient algebra turns out to be the right place to study finite
sections.  Note that if $T \in B(H)$ then $\sup \|P_n T P_n \| <
\infty$ (actually, $= \|T\|$) and hence the sequence of finite
sections $(P_n TP_n)$ can be regarded as a single element in $\Pi
M_{k(n)} ({\mathbb C})$.  The real element of interest, however, turns
out to be $\pi( (P_nTP_n))$.

\begin{center}
{\em Definitions arising from the finite section method.} 
\end{center}

\begin{defn}
Let $T \in B(H)$, $P_1 \leq P_2 \leq \cdots$ be a filtration of $H$
and $(P_n TP_n)$ be the corresponding finite sections.
\begin{enumerate}
\item $(P_n TP_n)$ is said to be {\em applicable} if there exists an
integer $n_0$ such that for every vector $y \in H$ and integer $n \geq
n_0$ there exists a unique solution $x_n$ to the equation $P_n T P_n
x_n = y_n$ and, moreover, the vectors $x_n$ converge (in norm) to a
solution $x$ of the equation $Tx = y$.  

\item $(P_n TP_n)$ is said to be {\em stable} if there exists an
integer $n_0$ such that for every integer $n \geq n_0$ the matrices
$P_n TP_n$ are invertible and $\sup_{n \geq n_0} \| (P_nTP_n)^{-1} \|
< \infty$.  A basic result of Polski states that $(P_n TP_n)$ is
applicable if and only if $T$ is invertible and $(P_n TP_n)$ is stable
(cf.\ \cite[Theorem 1.4]{HRS}).

\item The finite sections $(P_n T P_n)$ are said to be {\em fractal}
if for every subsequence $n_1 < n_2 < n_3 \cdots$ there exists a
$*$-isomorphism $\sigma : C^*(\pi((P_{n_j} T P_{n_j}))) \to
C^*(\pi((P_n T P_n)))$ such that $\sigma(\pi((P_{n_j} T P_{n_j}))) =
\pi((P_n T P_n))$. (We have abused notation here in using $\pi$ to
denote both of the canonical quotient maps $\Pi_n M_{k(n)} ({\mathbb
C}) \to \Pi_n M_{k(n)} ({\mathbb C})/ \oplus_n M_{k(n)} ({\mathbb C})$
and $\Pi_j M_{k(n_j)} ({\mathbb C}) \to \Pi_j M_{k(n_j)} ({\mathbb
C})/ \oplus_j M_{k(n_j)} ({\mathbb C})$.) 

\item For a matrix $S \in M_k ({\mathbb C})$ and $\varepsilon > 0$,
the {\em $\varepsilon$-regularization of $S$} is denoted
$S_{\varepsilon}$ and defined as follows.  Let $S = UDV$ be the
singular value decomposition of $S$ (i.e.\ $U,V \in M_k ({\mathbb C})$
are unitaries and $D = diag(\lambda_1, \ldots, \lambda_k)$ where
$\lambda_1 \leq \lambda_2 \leq \ldots \leq \lambda_k$ are the
(ordered) singular values of $S$.  Define $D_{\varepsilon} =
diag(\tilde{\lambda_1}, \ldots, \tilde{\lambda_k})$ where
$\tilde{\lambda_j} = \lambda_j$ if $\lambda_j > \varepsilon$ and
$\tilde{\lambda_j} = 0$ otherwise.  The $\varepsilon$-regularization
of $S$ is then defined to be $S_{\varepsilon} = UD_{\varepsilon}V$. 

\item The finite sections $(P_n TP_n)$ are said to be {\em stably
regularizable} if there exists a number $\varepsilon > 0$ such that
$\| (P_nTP_n)_{\varepsilon} - P_nTP_n\| \to 0$ and $\sup \|
(P_nTP_n)_{\varepsilon}^+ \| < \infty$.  In this case a Polski-type
theorem (cf.\ \cite[Theorem 2.13]{HRS}) states that for every vector
$y \in H$ the least square solutions $x_n$ to the (regularized)
equations $(P_nTP_n)_{\varepsilon} x_n = P_n y$ converge in norm to
the least square solutions of $Tx = y$ (i.e.\ $T$ is Moore-Penrose
invertible and $(P_nTP_n)_{\varepsilon}^+ y \to T^+ y$ in norm for
every $y \in H$).
\end{enumerate}
\end{defn}

One of the main topics of this paper is convergence aspects of
spectra.  As such it will be convenient to introduce some set notation
and recall Hausdorff convergence of sets.

\begin{center}
{\em Approximate inclusions and Hausdorff convergence of sets.} 
\end{center}

\begin{defn} Let ${\mathbb C}$ denote the complex plane.
\begin{enumerate}
\item For a set $\Sigma \subset {\mathbb C}$ and $\delta > 0$ we will
let $N_{\delta}(\Sigma)$ denote the $\delta$-neighborhood of $\Sigma$
(i.e.\ the union of all $\delta$-balls centered at points of $\Sigma$)
and $\overline{N_{\delta}(\Sigma)}$ be its closure. 

\item Given two sets $\Sigma, \Lambda \subset {\mathbb C}$ we say {\em
$\Sigma$ is $\delta$-contained in $\Lambda$} if $\Sigma \subset
N_{\delta}(\Lambda)$.  We prefer (and will use) the following notation
in this situation; $\Sigma \subset^{\delta} \Lambda$.

\item Given two compact sets $\Sigma, \Lambda \subset {\mathbb C}$
their {\em Hausdorff distance} is $$d_H (\Sigma, \Lambda) = \max \{
\sup_{\sigma \in \Sigma} d(\sigma, \Lambda), \sup_{\lambda \in
\Lambda} d(\lambda, \Sigma) \},$$ where $d(\sigma, \Lambda) =
\inf_{\lambda \in \Lambda} |\sigma - \lambda|$.  Note that $\delta >
d_H (\Sigma, \Lambda)$ if and only if $\Sigma$ is $\delta$-contained
in $\Lambda$ and vice versa -- i.e.\ $\Lambda \subset^{\delta} \Sigma$
and $\Sigma \subset^{\delta} \Lambda$.

\item Given a sequence of compact sets $\Sigma_n \subset {\mathbb C}$
and another compact $\Sigma \subset {\mathbb C}$ we will write
$$\Sigma_n \to \Sigma$$ if $d_H(\Sigma_n, \Sigma) \to 0$ -- i.e.\ if
for each $\delta > 0$ we have that $\Lambda_n \subset^{\delta}
\Lambda$ and $\Lambda \subset^{\delta} \Lambda_n$ for all sufficiently
large $n$.
\end{enumerate}
\end{defn}

We close this section by stating all of the theorems from \cite{HRS}
that we will need.  Roughly speaking these results say that all
convergence questions for the finite section method have a
C$^*$-formulation in $\Pi M_{k(n)} ({\mathbb C}) / \oplus M_{k(n)}
({\mathbb C})$.

\begin{thm}(Operator Equations)
\label{thm:OEHRS}
\begin{enumerate}
\item (Kozak's Theorem) $(P_n T P_n)$ is stable if and only if $\pi(
(P_n T P_n) )$ is an invertible element in $\Pi M_{k(n)} ({\mathbb C})
/ \oplus M_{k(n)} ({\mathbb C})$.

\item $(P_n T P_n)$ is stably regularizable if and only if $\pi( (P_n
T P_n) )$ is Moore-Penrose invertible in $\Pi M_{k(n)} ({\mathbb C}) /
\oplus M_{k(n)} ({\mathbb C})$.  
\end{enumerate}
\end{thm}

See \cite[Theorem 1.15]{HRS} for a proof of the first statement above
and \cite[Theorem 2.19]{HRS} for the second. 

\begin{thm}\cite[Theorem 2.28]{HRS} 
\label{thm:gcnHRS}
If $(P_nTP_n)$ is fractal and stably regularizable then the
generalized condition number of $\pi( (P_n T P_n) )$ is equal to $\lim
\| P_nTP_n \| \| X_n \|,$ where $(X_n)$ is any sequence such that
$\pi((X_n))$ is the Moore-Penrose inverse of $\pi( (P_n T P_n) )$.
\end{thm}

\begin{thm}(Spectral Approximations)
\label{thm:spectrumHRS} If $(P_nTP_n)$ is fractal then the following hold:
\begin{enumerate}
\item For every $\varepsilon > 0$, $\sigma^{(\varepsilon)}(P_n TP_n)
\to \sigma^{(\varepsilon)}(\pi((P_n TP_n))).$

\item $\sigma_2(P_n TP_n)
\to \sigma_2(\pi((P_n TP_n))).$

\item If each $P_n TP_n$ happens to be normal then $\sigma(P_n TP_n) \to
\sigma(\pi((P_n TP_n))).$ 
\end{enumerate}
\end{thm}

The first statement above is a combination of Theorems 3.31, 3.33(b)
and Proposition 3.6 in \cite{HRS}.  The second is a consequence of
Theorems 3.22, 3.23 and Proposition 3.6 in \cite{HRS}. The last
statement follows from Corollary 3.18, Theorem 3.20 and Proposition
3.6 in \cite{HRS}. 
%%%%%%%%%%%%%%%%%%%%%%%%%%%%%%%%%%%%%%%%%%%%%%%%%%%%%%%%%%%%%%%%%%%%%
\section{Existence and Uniqueness Results}

Let $P_1 \leq P_2 \leq \cdots$ be a filtration of $H$.  If the rank of
$P_n$ is some integer $k(n)$ and we identify $P_n B(H) P_n$ with the
complex $k(n) \times k(n)$ matrices $M_{k(n)} ({\mathbb C})$ then we
get a natural linear mapping $\phi_n : B(H) \to M_{k(n)} ({\mathbb
C})$ defined by $\phi_n(X) = P_n X P_n$ for all $X \in B(H)$.  We
further define linear maps $\Phi : B(H) \to \Pi M_{k(n)} ({\mathbb
C})$, $\Phi(X) = (\phi_n (X))$ and $\rho : B(H) \to \Pi M_{k(n)}
({\mathbb C})/ \oplus M_{k(n)} ({\mathbb C})$, $\rho = \pi \circ \Phi$
where $\pi : \Pi M_{k(n)} ({\mathbb C}) \to \Pi M_{k(n)} ({\mathbb
C})/ \oplus M_{k(n)} ({\mathbb C})$ is the canonical quotient mapping.

The following lemma, which is well known to the operator algebra
community, is really the key observation of this section.  The proof
is also very simple so we include it for the reader's convenience.

\begin{lem}  
\label{thm:QDremark}
If $\| [T, P_n] \| \to 0$ then $\rho$
restricts to an injective $*$-homomorphism on the unital
C$^*$-algebra, denoted C$^*(T)$, generated by $T$. In other words,
$\rho|_{C^*(T)} : C^*(T) \to \Pi M_{k(n)} ({\mathbb C})/ \oplus
M_{k(n)} ({\mathbb C})$ is a $*$-isomorphism onto it's image.
\end{lem}

\begin{proof} First note that the projections $P_n$ asymptotically
commute with $T^*$ as well (since $\| [T, P_n] \| = \| [T, P_n]^* \| =
\| [T^*, P_n] \|$).  A simple interpolation argument together with the
triangle inequality shows that $(P_n)$ asymptotically commute with
any polynomial in $T$ and $T^*$.  It follows that the $P_n$'s
asymptotically commute with every element in $C^*(T)$.  Since $\|
\rho(X) \| = \| X \|$ for all $X \in B(H)$ (in particular, $\rho$ is
injective) it only remains to show that $\rho$ restricts to a
multiplicative map on $C^*(T)$.  This follows easily since for all $a,
b \in C^*(T)$ we have $\| \phi_n(ab) - \phi_n(a)\phi_n(b) \| = \| P_n
ab P_n - P_n a P_n b P_n \| = \| P_n (P_n a - a P_n) bP_n \| \to 0$.
\end{proof}

\begin{lem}
\label{thm:fractal}
Let $T \in B(H)$ be a quasidiagonal operator and $P_1 \leq P_2
\leq \cdots$ be any increasing sequence of finite rank projections
such that $\| [P_n, T] \| \to 0$ and $\| P_n(v) - v \| \to 0$, for all
$v \in H$, as $n \to \infty$. Then $(P_n T P_n)$ is fractal.
\end{lem}

\begin{proof} 
Every subsequence of $\{ P_n \}$ also asymptotically commutes with $T$
and converges strongly to the identity.  Hence the corresponding map
$\rho : C^*(T) \to \Pi M_{k(n_j)} ({\mathbb C})/ \oplus M_{k(n_j)}
({\mathbb C})$ is a $*$-isomorphism.
\end{proof}

As stated in the title of this section, we regard the following results 
as  existence and uniqueness results.  That is, for a quasidiagonal
operator one can always find a sequence of projections such that the
finite sections converge where we want (existence) and, in fact, the
finite sections arising from  any other {\em asymptotically commuting}
sequence will converge to the same thing (uniqueness).

\begin{thm}
\label{thm:existOE}
(Operator Equations) Let $T \in B(H)$ be a quasidiagonal operator and
$P_1 \leq P_2 \leq \cdots$ be any filtration of $H$ which
asymptotically commutes with $T$. Then the following statements hold:
\begin{enumerate}
\item The sequence $(P_n T P_n) \in \Pi M_{k(n)} ({\mathbb C})$ is
stable if and only if $T$ is
invertible. 

\item The sequence $(P_n T P_n) \in \Pi M_{k(n)}
({\mathbb C})$ is applicable if and only if $T$ is invertible.

\item The sequence $(P_n T P_n)$ is stably regularizable if and only
if $T$ is Moore-Penrose invertible (i.e.\ the range of $T$ is closed
-- see \cite[Theorem 2.4]{HRS}).

\item
Assume that $T$ has closed range, $\varepsilon > 0$ is some number
strictly smaller than $\inf \{ \lambda : 0 < \lambda \in
\sigma_2(T)\}$, $y \in H$ is given and $x \in H$ is the least square
solution to the equation $Tx = y$.  If $y_n = P_n y$ and $x_n$ is the
least square solution to (the $\varepsilon$-regularized equation)
$(P_nTP_n)_{\varepsilon} x_n = y_n$ then $\| x_n - x \| \to 0$ as $n
\to \infty$.
\end{enumerate}
\end{thm}

\begin{proof} 
(1) follows from Kozak's Theorem and Lemma \ref{thm:QDremark} since
$\rho(T)$ is invertible if and only if $T$ is invertible.

(2) follows immediately from (1) and Polski's Theorem which
was stated right after the definition of stability in the previous
section.

(3) follows from the second part of Theorem \ref{thm:OEHRS} and Lemma
    \ref{thm:QDremark}.

(4) follows from (3) and the least square version of
Polski's Theorem stated after the definition of stably regularizable.
Strictly speaking, \cite[Theorem 2.13]{HRS} only ensures the existence
of an $\varepsilon > 0$ for which our claim holds true.  However, the
convergence of spectra results below will ensure that any $0 <
\varepsilon < \inf \{ \lambda : 0 < \lambda \in \sigma(T^* T)\}$ will
do.
\end{proof}

\begin{thm}(Norms and Condition Numbers)
\label{thm:existnorms}
Let $T \in B(H)$ be a quasidiagonal operator and
$P_1 \leq P_2 \leq \cdots$ be any filtration of $H$ which
asymptotically commutes with $T$. Then the following statements hold:
\begin{enumerate}
\item $\| T \| = \lim\limits_{n \to \infty} \| P_n T P_n \|$ and $\|
T^{-1} \| = \lim\limits_{n \to \infty} \| (P_n T P_n)^{-1} \|$.

\item If $T$ is invertible then $cond(T) = \lim\limits_{n \to \infty}
cond(P_n T P_n)$. 

\item If $T$ has closed range, $\varepsilon > 0$ is some number
strictly smaller than $\inf \{ \lambda : 0 < \lambda \in
\sigma_2(T)\}$, and $T^{+}$ (resp.\ $(P_nTP_n)_{\varepsilon}^{+}$)
denotes the Moore-Penrose inverse of $T$ (resp.\
$(P_nTP_n)_{\varepsilon}$) then the generalized condition numbers of
$(P_nTP_n)_{\varepsilon}$ converge to the generalized condition
numbers of $T$ -- i.e.\ $\| T \| \| T^{+} \| = \lim \|
(P_nTP_n)_{\varepsilon} \| \| (P_nTP_n)_{\varepsilon}^{+} \|$.
\end{enumerate}
\end{thm}

\begin{proof}
The first part of statement (1) is trivial. In fact, for every $X
\in B(H)$, the norms of $P_n X P_n$ {\em increase} up to $\| X \|$
(since $P_n P_{n+1} = P_n$ and thus $\| P_n X P_n \| = \| P_n
(P_{n+1} X P_{n+1}) P_n \| \leq \|P_{n+1} X P_{n+1}\|$).  For the
inverse statement we first assume that $T$ is invertible and note that
$\| T^{-1} \| = \lim \| P_n T^{-1} P_n \|$.  Also, since $T^{-1} \in
C^*(T)$, the fact that $\rho(T^{-1}) = \rho(T)^{-1}$ implies that
$(P_n T^{-1} P_n) - ((P_n T^{-1} P_n)^{-1}) \in \oplus M_{k(n)}
({\mathbb C})$. It follows that $\| P_n T^{-1} P_n - (P_n T^{-1} P_n)^{-1}
\| \to 0$ as $n \to \infty$ and hence $\| T^{-1} \| = \lim \| P_n
T^{-1} P_n \| = \lim \| (P_n T P_n)^{-1} \|$.  If $T$ is not
invertible but $\lim\limits_{n \to \infty} \| (P_n T P_n)^{-1} \| \neq
\infty$ then we can find a subsequence $\{ n_j \}$ such that
$\lim\limits_{j \to \infty} \| (P_{n_j} T P_{n_j})^{-1} \| = c <
\infty$.  However, the projections $\{ P_{n_j} \}$ also
asymptotically commute with $T$ and the finite sections
$(P_{n_j}TP_{n_j})$ are now stable.  Thus, by Theorem
\ref{thm:existOE} (1), $T$ is invertible -- a contradiction.

(2) follows immediately from (1). 

(3) follows from Theorem \ref{thm:gcnHRS} and Lemmas
    \ref{thm:QDremark} and \ref{thm:fractal}.
\end{proof}

\begin{thm}(Spectral Approximations) 
\label{thm:existspectrum}
Let $T \in B(H)$ be a quasidiagonal operator and
$P_1 \leq P_2 \leq \cdots$ be any filtration of $H$ which
asymptotically commutes with $T$. Then the following statements hold:
\begin{enumerate}
\item For every $\varepsilon > 0$, $\sigma^{(\varepsilon)}(P_n T P_n)
\to \sigma^{(\varepsilon)}(T)$.

\item $\sigma(T) = \bigcap\limits_{\varepsilon > 0} \bigg( \lim
\sigma^{(\varepsilon)}(P_n T P_n)\bigg),$ where $\lim
\sigma^{(\varepsilon)}(P_n T P_n) = \sigma^{(\varepsilon)}(T)$ by (1). 

\item $\sigma_2(P_n T P_n) \to \sigma_2 (T)$.
\end{enumerate}
\end{thm}

\begin{proof} Given Lemmas \ref{thm:QDremark} and \ref{thm:fractal},
the proof is reduced to recalling Theorem \ref{thm:spectrumHRS} from
the last section.
\end{proof}

\begin{rem} As mentioned earlier, a theorem of Berg states that every
normal operator on a Hilbert space is quasidiagonal and hence the
existence and uniqueness results above can be applied.  Even in the
case of a self-adjoint operator the {\em existence} of a filtration for
which the finite sections are guaranteed to converge nicely appears to
be a new result. 
\end{rem}

\begin{rem} 
Many standard examples of operators are not quasidiagonal.  For
example, any Fredholm operator with non-zero index is not
quasidiagonal.  However, there are very nice characterizations of
quasidiagonality for essentially normal operators (cf.\
\cite[Corollary IX.7.4, Section IX.8]{davidson}) and weighted shift
operators (\cite{smucker}, \cite{narayan}).
\end{rem}

\begin{rem} 
Though many operators are not quasidiagonal, a remarkable result of
Voiculescu implies that if one starts with an arbitrary $T \in B(H)$
then it is possible to construct a quasidiagonal operator out of $T$.
Namely, let $\lambda_1, \lambda_2, \ldots$ be any listing of all the
rational numbers in the half open interval $(0,1]$.  If $T \in B(H)$
is arbitrary then Voiculescu's homotopy invariance theorem (cf.\
\cite{dvv:QDsurvey}) implies that
$$\tilde{T} = \oplus_{i \in {\mathbb N}} \lambda_i T,$$ acting on
$\oplus_{i \in {\mathbb N}} H$, is a quasidiagonal operator. 
\end{rem}

We close this section with a result which has nothing to do with
quasidiagonality.  Namely, we will consider the convergence of the
algebraic numerical ranges of finite sections.  It is well known that
spatial numerical ranges are well behaved (cf.\ \cite[Theorem
3.52]{HRS}), but algebraic numerical ranges also behave nicely.
Though we would be a bit surprised if this result is not already
known, it is not mentioned in \cite{HRS} and hence we include a proof.

\begin{thm}(Numerical Ranges) 
Let $T \in B(H)$ be arbitrary and $P_1 \leq P_2 \leq
\ldots$ be any sequence of finite rank projections such that $\|
P_n(v) - v \| \to 0$, as $n \to \infty$, for all $v \in H$.  Then the
following statements hold:
\begin{enumerate}
\item For every $n$, $SN(P_nTP_n) = AN(P_nTP_n) \subset SN(T) \subset AN(T)$.

\item $SN(P_nTP_n) = AN(P_n T P_n) \to AN(T) = \overline{SN(T)} = AN(T)$.
\end{enumerate}
\end{thm}

\begin{proof} 
It is immediate from the definitions that $SN(X) \subset AN(X)$ for
any $X \in B(H)$. Moreover, it is also known that $AN(X) =
\overline{SN(X)}$ (cf.\ \cite[Theorem 3.45]{HRS}).  Since the unit
ball of a finite dimensional Hilbert space (e.g.\ $P_nH$) is (norm)
compact, these two facts imply that $SN(P_nTP_n) = AN(P_nTP_n)$ (for a
given state $\phi$, choose unit vectors $v_k \in P_nH$ such that
$\phi(P_nTP_n) = \lim \langle P_nTP_n v_k,v_k \rangle$ and any cluster
point of $\{v_k\}$ will do the trick).  The proof of (1) is then
complete after we observe that $SN(P_nTP_n) \subset SN(T)$: for a unit
vector $v \in P_nH$ we have $\langle Tv, v \rangle = \langle TP_nv,
P_nv \rangle = \langle P_nTP_nv, v \rangle$.

Statement (2) follows from (1) and the general fact that $SN(P_n T
P_n) \to \overline{SN(T)}$ (cf.\ \cite[Theorem 3.52]{HRS}).
\end{proof}

%%%%%%%%%%%%%%%%%%%%%%%%%%%%%%%%%%%%%%%%%%%%%%%%%%%%%%%%%%%%%%%%%%
\section{Commutators and Rates of Convergence}

In this section we observe that norms of the commutators $\| [P_n,
T]\|$ provide a simple and natural upper bound for the rate of
convergence of finite sections $(P_n T P_n)$.  For spectral
approximations this only yields `one sided' rates of convergence but,
sadly, a trivial example at the end of this section shows that this is
best possible.  

From a computational point of view this will, in general, not provide
a very useful measurement. However, we will see in later sections that
in some cases these commutator estimates yield computationally
relevant rates of convergence.

As was the case for Lemma \ref{thm:QDremark}, the following (well
known) lemma is really the key remark of this section.

\begin{lem} 
\label{thm:commremark}
Let $T \in B(H)$ be arbitrary and $P \in B(H)$ be a self-adjoint
projection.  Then $$\| T - (PTP + P^{\perp}TP^{\perp}) \| \leq
\|PT - TP\|,$$ where $P^{\perp} = I - P$.
\end{lem}
 
\begin{proof} Since $T = PTP + PTP^{\perp} + P^{\perp}TP +
P^{\perp}TP^{\perp}$ ($T = I T I$ and $I = P + P^{\perp}$) we have
that $T - (PTP + P^{\perp}TP^{\perp}) = PTP^{\perp} + P^{\perp}TP$.
Since the domains and ranges of the operators $PTP^{\perp}$ and
$P^{\perp}TP$ are orthogonal we have $\|PTP^{\perp} +
P^{\perp}TP\| = \max \{ \|PTP^{\perp}\|, \|P^{\perp}TP\| \}$.
However, each of the latter norms is dominated by $\|PT - TP\|$ -- for
example, $\| P^{\perp}TP \| = \| P^{\perp}TP - P^{\perp}PT \| \leq \|
P^{\perp} \| \| TP - PT\|$.
\end{proof}

Though we prefer to state the results below in terms of norms of
commutators, we should also point out that these norms are determined
by the norms of finite dimensional matrices:
\begin{eqnarray*}
\| PT - TP \|  & = & \| PTP^{\perp} - P^{\perp}TP \|\\ 
& = &  \max\{ \| PTP^{\perp} \|, \| P^{\perp}TP \| \}\\
& = &  \max \{ \| P^{\perp}T^*P \|, \| P^{\perp}TP \| \}\\
& = &  \max \{ \|PT P^{\perp}T^*P \|^{1/2}, \| PT^*P^{\perp}TP \|^{1/2} \}.
\end{eqnarray*}

\begin{prop}(Operator Equations) 
\label{thm:rateOE}
Let $T \in B(H)$ be a quasidiagonal operator and $P_1 \leq P_2
\leq \cdots$ be any filtration 
such that $\| [P_n, T] \| \to 0$. Then the following statements hold:
\begin{enumerate}
\item If $T$ is invertible and $\| [P_n, T] \| < \frac{1}{\| T^{-1}
\|}$ then $P_n T P_n$ is also invertible. 

\item If $T$ is invertible, $\| [P_n, T] \| < \frac{1}{\| T^{-1} \|}$,
$y \in H$ is given, $x \in H$ is the unique solution to $Tx = y$, $y_n
= P_n y$ and $x_n$ is the unique solution to $P_n T P_n x_n = y_n$
then $$\| x - x_n \| \leq \| T^{-1} \| \big( \| y - y_n \| + \| x_n \|
\| [P_n, T] \| \big).$$
\end{enumerate}
\end{prop}

\begin{proof} The first statement follows from the following
inequalities: $$\| P_nTP_n T^{-1} P_n - P_n \| = \| P_n (TP_n - P_n
T)T^{-1}P_n \| \leq \| [P_n, T]\| \| T^{-1} \| < 1.$$ Indeed, this
shows that $P_nTP_n T^{-1} P_n = (P_nTP_n)(P_nT^{-1}P_n)$ is
invertible (we now regard $P_n$ as the unit of the matrix algebra
$P_nB(H)P_n$) and hence $P_nTP_n$ is right invertible.  However, in
finite dimensions right invertible implies invertible.

The second assertion is also straightforward: 
\begin{eqnarray*}
\| x_n - x\|  & = & \| T^{-1}(Tx_n - Tx)\|\\ 
& \leq &  \|T^{-1}\| \|Tx_n - y \|\\
& = &  \|T^{-1}\|\|TP_nx_n - P_nTP_nx_n + y_n - y\|\\
& \leq & \|T^{-1}\|\big( \| y_n - y\| + \| x_n \| \| 
(TP_n - P_nT)P_n \|\big)\\
& \leq & \| T^{-1} \| \big( \| y - y_n \| + \| x_n \|
\| [P_n, T] \| \big).
\end{eqnarray*}
\end{proof}

\begin{rem}({\em Computability}) We regard the estimates above as
reasonably computable objects (at least in concrete cases).  Indeed,
if $T$ and $\{ P_n \}$ are concretely given then $\| [P_n, T] \|$ is
computable in terms of two finite dimensional matrices. (Warning: In
this remark we take the point of view that anything in finite
dimensions is `computable' and limits of finite dimensional
calculations are regarded as `reasonably computable'.)  Also, note 
that $\| y - P_n y \| = \sqrt{\|y\|^2 - \|P_n y\|^2}$
is computable (at least if one knows $\|y\|$).

Recall that the polar decomposition implies that for any invertible
operator $S$ we have $$\| S^{-1} \| = \sup_{\lambda \in \sigma_2(S)}
\frac{1}{\lambda}.$$ Hence the quantities $\| T^{-1} \|$ in the
estimates above are reasonably computable since we know the singular
values of $P_nTP_n$ converge to the singular values of $T$.  Also, $\|
x_n \| \leq \| (P_nTP_n)^{-1} \| \| y_n \|$ and thus can be estimated
in concrete cases.
\end{rem}

\begin{rem}({\em Least Square Estimates}) We have been unable to find
any `reasonable' estimates on the rate of convergence of least square
solutions.  Our feeling is that this may be possible, however we don't
yet see how to avoid the Moore-Penrose projection (a certain spectral
projection for $T^* T$ -- see \cite[Theorem 2.15]{HRS}) which arises
in the construction of $T^+$. 
\end{rem}

\begin{prop}(Spectral Approximations) 
\label{thm:ratespectrum}
Let $T \in B(H)$ be a quasidiagonal operator and $P_1 \leq P_2
\leq \cdots$ be any filtration 
such that $\| [P_n, T] \| \to 0$. Then the following statements hold:
\begin{enumerate}
\item For every $\varepsilon > 0$, $\sigma^{(\varepsilon)}(P_n TP_n)
\subset \sigma^{(\varepsilon + \|[P_n, T]\|)}(T)$.

\item $\sigma_2 (P_nTP_n) \subset^{\delta_n} \sigma_2(T)$, where
$\delta_n = \sqrt{2\|T\|\|[P_n,T]\|}$.
\end{enumerate}
\end{prop}

\begin{proof} We will need the following characterization of
$\varepsilon$-pseudospectra (cf.\ \cite[Theorem 3.27]{HRS}): If $X \in
B(H)$ then $$\sigma^{(\varepsilon)} (X) = \bigcup_{Y \in B(H), \ \|Y\|
\leq \varepsilon} \sigma(X + Y).$$ Given this fact, it is clear that
if $\| X - Z \| \leq \delta$ then $\sigma^{(\varepsilon)} (X) \subset
\sigma^{(\varepsilon + \delta)} (Z)$.  Since $\sigma^{(\varepsilon)}
(P_nTP_n) \subset \sigma^{(\varepsilon)} (P_nTP_n +
P_n^{\perp}TP_n^{\perp})$, statement (1) follows from Lemma
\ref{thm:commremark} and the remarks above. 

Part (2) requires the following well known fact: If $X, Y \in B(H)$,
$X$ is self-adjoint (even normal will suffice) and $\| X - Y \| \leq
\delta$ then $\sigma(Y) \subset^{\delta} \sigma(X)$. To complete the
proof we first note that $\sigma(P_nT^*P_nTP_n) \subset
\sigma(P_nT^*P_nTP_n + P_n^{\perp}T^*P_n^{\perp}TP_n^{\perp}) =
\sigma((P_nT^*P_n + P_n^{\perp}T^*P_n^{\perp})(P_nTP_n +
P_n^{\perp}TP_n^{\perp}))$.  However, Lemma \ref{thm:commremark}
implies that $\| (P_nT^*P_n + P_n^{\perp}T^*P_n^{\perp})(P_nTP_n +
P_n^{\perp}TP_n^{\perp}) - T^*T \| \leq 2\|T\|\|[P_n, T]\|$.
Combining all these facts we see that $\sigma(P_nT^*P_nTP_n)
\subset^{2\|T\|\|[P_n, T]\|} \sigma(T^*T)$ -- i.e. if $\lambda \in
\sigma_2(P_nTP_n)$ then there exists some point $\gamma \in \sigma_2
(T)$ such that $| \lambda^2 - \gamma^2 | \leq 2\|T\|\|[P_n, T]\|$.
Since $(\lambda - \gamma)^2 \leq | \lambda^2 - \gamma^2 |$, statement
(2) follows.  Note that for large singular values one does not have to
take square roots, however the estimate $(\lambda - \gamma)^2 \leq |
\lambda^2 - \gamma^2 |$ becomes sharp at zero.
\end{proof}

We end this section with a simple example which shows that there is no
hope of giving rates of convergence of the quantities not mentioned in
this section in terms of norms of commutators. 

Let $j_1 < j_2 < j_3 \cdots$ be integers and consider the diagonal,
self-adjoint operator on $l^2 ({\mathbb N})$ whose matrix (w.r.t.\ the
canonical orthonormal basis) is $$diag(2,\ldots,2,1 + 1/2,3 -
1/2,2,\ldots,2,1+1/3,3-1/3,2,\ldots,2,1+1/4,3-1/4,2\cdots),$$ where
the first string of 2's has length $j_1$, the second string of 2's has
length $j_2$ and so on.  One should have in mind here that the $j_k$'s
are growing like exponentials (or faster, if you like) in $k$.  Now
let $P_n$ denote the orthogonal projection onto the span of the first
$n$ basis vectors of $l^2 ({\mathbb N})$.  

Note that $[P_n, T] = 0$ for all $n$.  On the other hand, the norms of
$P_n T P_n$ and it's inverse approach the norms of $T$ and it's
inverse as slowly (in $n$) as you like.  Also, the condition numbers
of $P_n TP_n$ are approaching the condition number of $T$ as slowly as
you like. 

It would also be nice if we could approximately reverse the inclusions
in Proposition \ref{thm:ratespectrum} -- i.e.\ if we knew something
like $\sigma(T) \subset \sigma^{(\varepsilon)}(P_nTP_n)$ for all $n$
larger than some number which depended on commutators.  This would
amount to saying how quickly the pseudo-spectra of $P_nTP_n$ ``fill
up'' the pseudo-spectrum of $T$.  Unfortunately, the example above
shows that this can happen as slowly as desired.  Similarly the
numerical ranges of $P_nTP_n$ will ``fill up'' the numerical range of
$T$ at as slow a rate as one wishes to prescribe.

Finally, it would be nice for solving operator equations if we could
somehow relate $\| y - y_n \|$ and $\| [P_n, T] \|$ but it is not hard
to see that this is also impossible to do in general.

%%%%%%%%%%%%%%%%%%%%%%%%%%%%%%%%%%%%%%%%%%%%%%%%%%%%%%%%%%%%%%%%%%
\section{Normal and Self-Adjoint Operators}

Here we observe that some of the statements from the previous section
can be cleaned up a bit in the case that $T$ is a normal or, better
yet, self-adjoint operator.  At the end of this section we also state
a result from \cite{brown:AFD} which provides a very general {\em
existence} result which is analogous to the Szeg\"{o}-type results
obtained by Bill Arveson for certain self-adjoints (cf.\ \cite[Chapter
7]{HRS}).

\begin{prop} If $T \in B(H)$ is a normal operator and $P_1 \leq P_2
\leq \cdots$ is any filtration such that $\| [P_n, T] \| \to 0$ then:
\begin{enumerate}
\item For every $\varepsilon > 0$, $\sigma^{(\varepsilon)}(P_n T P_n)
\to \overline{N_{\epsilon}(\sigma(T))}$.

\item For every $\varepsilon > 0$, $\sigma^{(\varepsilon)}(P_n TP_n)
\subset^{\delta_n} \sigma(T)$, where $\delta_n =
\varepsilon + \| [T, P_n]\|$.
\end{enumerate}
\end{prop}

\begin{proof} For a normal operator $T$ we have
$\overline{N_{\epsilon}(\sigma(T))} = \sigma^{(\varepsilon)}( T )$.
This follows from the fact that if $T$ is normal and $\lambda \in
{\mathbb C}$ then $\| (\lambda - T)^{-1} \| = \frac{1}{d(\lambda,
\sigma(T))}$, where $d(\lambda, \sigma(T))$ is the distance from
$\lambda$ to $\sigma(T)$.  Hence we have just reformulated the results
of sections 2 and 3.
\end{proof}

\begin{rem}  
On page 116 in \cite{HRS} a general conjecture concerning convergence
of spectra of `essentially normal' approximation sequences is
formulated.  Though we have been unable to resolve this conjecture we
do wish to point out that an affirmative answer would lead to a nice
improvement in the previous result.  Namely, it would follow that
$\sigma(P_n T P_n) \to \sigma(T)$ -- i.e.\ we wouldn't have to bother
with pseudospectra in the result above.
\end{rem}

\begin{prop} If $T \in B(H)$ is a self-adjoint operator and $P_1 \leq P_2
\leq \cdots$ is any filtration such that $\| [P_n, T] \| \to 0$ then:
\begin{enumerate}
\item $\sigma(P_n T P_n) \to \sigma(T)$.

\item $\sigma(P_n TP_n) \subset^{\delta_n} \sigma(T)$,
where $\delta_n = \| [T, P_n]\|$.
\end{enumerate}
\end{prop}

\begin{proof} Assertion (1) follows from part (3) of Theorem
\ref{thm:spectrumHRS} and Lemmas \ref{thm:QDremark} and
\ref{thm:fractal}.  The second part follows from Lemma
\ref{thm:commremark}, the fact that $\sigma(P_n TP_n) \subset
\sigma(P_n TP_n + P^{\perp}_n TP^{\perp}_n)$ and the first sentence in
the proof of part (2) of Proposition \ref{thm:ratespectrum}.
\end{proof}

\begin{rem} Following Arveson's ideas, it is possible to also compute
the essential spectrum of a self-adjoint operator via quasidiagonal
finite sections (cf.\ \cite[Chapter 7]{HRS}).  It does not seem easy,
however, to develop any rate of convergence estimates in this case. 
\end{rem}

We conclude this section with an {\em existence} result which is
analogous to Szeg\"{o}'s celebrated limit theorem concerning spectral
distributions of self-adjoint Toeplitz operators.  In
\cite{brown:AFD} we studied various notions of invariant means
on C$^*$-algebras and observed in Theorem 10.4 that `quasidiagonal'
invariant means can be recovered following the ideas first laid out by
Bill Arveson (as explained in \cite[Section 7.2]{HRS}).  To avoid
technicalities we will just restate the result in the case of normal
operators. 

\begin{thm} Let $T \in B(H)$ be a normal operator.  Then there exists
a filtration $P_1 \leq P_2 \leq \cdots$ such that $\| [P_n, T]\| \to
0$ and with the property that for every Borel, probability measure
$\mu$ supported on the essential spectrum of $T$ there exists a
subsequence $\{ P_{n_k} \}$ such that $$tr(P_{n_k} \cdot P_{n_k}) \to
\mu$$ in the weak$-*$ topology, where $tr(\cdot)$ is the (unique)
tracial state on $P_{n_k}B(H)P_{n_k}$ and $tr(P_{n_k} \cdot P_{n_k})$
is the state on $C^*(T)$ defined by $x \mapsto tr(P_{n_k} x P_{n_k})$.
\end{thm}

\begin{proof} Since normal operators generate exact C$^*$-algebras
(even nuclear ones) and every trace on an abelian C$^*$-algebra is a
quasidiagonal invariant mean (cf.\ \cite[Lemma 6.4]{brown:AFD})
this result is a special case of \cite[Theorem 10.4]{brown:AFD}.
\end{proof}

%%%%%%%%%%%%%%%%%%%%%%%%%%%%%%%%%%%%%%%%%%%%%%%%%%%%%%%%%%%%%%%%%%
\section{Improved Existence Result when $C^*(T)$ is Exact}

Mainly due to the pioneering work of Kirchberg, {\em exact}
C$^*$-algebras have been at the forefront of operator algebra theory
for the past decade.  (See \cite{kirchberg}, \cite{wassermann}.)  This is a
technical notion which we do not wish to discuss here, however it
turns out that most (but not all) of the basic examples of
C$^*$-algebras enjoy this property (commutative algebras, AF algebras,
irrational rotation algebras, Cuntz algebras and all subalgebras of
these examples).  Similarly, it turns out that many standard examples
of operators on Hilbert space generate C$^*$-algebras with this
property.  For example, $C^*(T)$ is exact whenever $T \in B(H)$ is
essentially normal (e.g.\ self-adjoint, unitary or normal) or if $T$
is a band operator (either unilateral or bilateral).  Moreover, it
follows from Arveson's extension theorem that this class of operators
is closed under norm limits and hence band dominated operators (i.e.\
norm limits of band operators) also generate exact C$^*$-algebras.

In the early 80's Herrero asked whether every quasidiagonal operator
$T$ was the norm limit of operators $T_n$ such that $C^*(T_n)$ was
finite dimensional for each $n$.  This turns out to be false (cf.\
\cite{szarek}, \cite{voiculescu}) but the precise description of
operators with this approximation property was discovered a couple
years ago (cf.\ \cite{brown:Herrero}, \cite{dadarlat}).

\begin{thm} Let $T \in B(H)$ be given.  There exist operators
$T_n \in B(H)$ such that $\| T - T_n \| \to 0$ and $C^*(T_n)$ is
finite dimensional (for all $n$) if and only if $T$ is quasidiagonal
and $C^*(T)$ is exact. 
\end{thm}

This result is relevant to this paper because of the simple structure
of finite dimensional C$^*$-algebras: they are just finite direct sums
of finite dimensional matrix algebras.  So if $T$ is close to $T_n$
and $C^*(T_n)$ is finite dimensional then one immediately gets a
matrix whose $\varepsilon$-pseudospectrum, for example, approximates
that of $T$ (since we may regard $T_n$ as a direct sum of finite
dimensional matrices when viewed as an element of $C^*(T_n)$).  In
this setting rates of convergence are controlled by $\| T_n - T \|$
and estimates on the Hausdorff distance of spectra (in the normal
case, at least) can be given.  Moreover, standard representation
theory of finite dimensional C$^*$-algebras gives a decomposition of
the underlying Hilbert space as a direct sum of finite dimensional
reducing subspaces for $C^*(T_n)$ and thus a way of getting
appropriate finite sections.  Of course, actually finding these
approximating $T_n$'s and corresponding finite sections in specific
examples does not appear to be easy and hence we have little more to
say about it at this point.  However, an effective means of computing
spectra, inspired by the observations of this section, in irrational
rotation algebras is worked out in \cite{brown:PV}.

%%%%%%%%%%%%%%%%%%%%%%%%%%%%%%%%%%%%%%%%%%%%%%%%%%%%%%%%%%%%%%%%%%
\section{Quasidiagonal  Unilateral Band Operators}

Though the results of the previous sections are quite satisfactory
from a theoretical point of view some practical numerical analytic
questions have went unaddressed.  For example, if an explicitly given
operator is known to be quasidiagonal how does one explicitly write
down asymptotically commuting finite rank projections?  Can one
control the ranks of these projections?  In the generality of the
previous sections these questions cannot be answered
(even for a general self-adjoint operator this is probably hopeless).
However, in special cases some progress can be made and this is the
subject of the final two sections of this paper. 

In this section we consider certain band operators.  As mentioned
before, not all band operators are quasidiagonal (e.g.\ unilateral
shift) but a fairly broad class of band operators, described below,
can be attacked with this technology. 

Our Hilbert space will now be $l^2({\mathbb N})$ with canonical
orthonormal basis $\{ \delta_i \}_{i \in {\mathbb N}}$.  We will
regard $l^{\infty}({\mathbb N})$ as acting on $l^2({\mathbb N})$ by
multiplication operators (i.e.\ $f \in l^{\infty}({\mathbb N})$ acts
on $l^2({\mathbb N})$ by sending $\delta_i$ to $f(i)\delta_i$).  $T
\in B(l^2({\mathbb N}))$ will be a weighted shift operator with weight
sequence $(\alpha_1, \alpha_2, \ldots)$ -- i.e.\ $(\alpha_1, \alpha_2,
\ldots)$ is a bounded sequence of complex numbers and $T(\delta_i) =
\alpha_i \delta_{i+1}$.  Our basic assumption is: 
\begin{center} {\em Assume the existence of a subsequence $\{ i_k\}_{k \in
{\mathbb N}}$ such that $|\alpha_{i_k}| \to 0$ as $k \to \infty$.}
\end{center} 
Note that a unilateral weighted shift is quasidiagonal if and only if
the assumption above holds.  Below we will observe that this
assumption implies quasidiagonality while the absence of such a
subsequence implies that the weighted shift is a Fredholm operator
with index $-1$ and, hence, is not quasidiagonal.

The band operators considered will be of the form $$S = f_{-m}(T^*)^m
+ f_{-m+1}(T^*)^{m-1} + \cdots f_{-1} T^* + f_0 + f_1 T + \cdots +
f_{m-1}T^{m-1} + f_mT^m,$$ where $f_j \in l^{\infty}({\mathbb N})$,
$-m \leq j \leq m$, and $m$ is a positive integer.  Note that $S$ is a
band operator whose main diagonal is just $f_0$.  The off diagonal
entries of $S$ are somewhat complicated to describe, but the
`coefficients' $f_j$ are arbitrary and hence operators of the form
above are still are rather large class.  The reader may find it useful
to consider a weight sequence of only zeroes and ones to get a feel
for the kinds of operators arising this way.  In fact, any band
operator whose matrix has ``butterflies'' tending to zero can be
handled with the quasidiagonal approach.  In other words, any band
operator which is a compact perturbation of a block diagonal band
operator can be handled this way but rates of convergence can be
complicated to state precisely so we stick to the situation described
above. 

For ease of notation we will let $q(x) = x^{m} + x^{m-1} + \cdots + x$
and $q^{\prime}$ be the derivative of this polynomial.  We will also
let $M = \max_{-m \leq j \leq m}\{ \| f_j \| \}$. Finally, we define
finite rank projections $P_k$ as the orthogonal projections onto the
span of $\{ \delta_1, \ldots, \delta_{i_k}\}$.

\begin{thm} 
\label{thm:unilateral}
Let $T$, $S$ and $\{ P_k \}$ be as above. Then the following
statements hold:
\begin{enumerate}
\item If $S$ is invertible and $|\alpha_{i_k}| <
\frac{1}{\|S^{-1}\|Mq^{\prime}(\|T\|)}$ then $P_kSP_k$ is invertible.

\item If $S$ is invertible, $|\alpha_{i_k}| <
\frac{1}{\|S^{-1}\|Mq^{\prime}(\|T\|)}$, $y \in l^2({\mathbb N})$ is
given, $x \in l^2({\mathbb N})$ is the unique solution to the equation
$Sx = y$, $y_k = P_k y$ and $x_k$ is the unique solution to $P_kSP_k
x_k = y_k$ then $$\| x_k - x \| \leq \| S^{-1} \| \big( \| y - y_n \|
+ \| x_n \|M q^{\prime}(\|T\|)|\alpha_{i_k}| \big).$$

\item For every $\varepsilon > 0$, the $\varepsilon$-pseudospectra of
$P_k S P_k$ converge to the $\varepsilon$-pseudospectrum of $S$ and
$\sigma^{(\varepsilon)}(P_n SP_n) \subset \sigma^{(\varepsilon + M
q^{\prime}(\|T\|)|\alpha_{i_k}|)}(S)$.

\item The singular values of $P_k SP_k$ converge to the singular
values of $S$ and $\sigma_2 (P_nTP_n) \subset^{\delta_n} \sigma_2(T)$,
where $\delta_n = \sqrt{2\|T\|M q^{\prime}(\|T\|)|\alpha_{i_k}|}$.
\end{enumerate}
\end{thm}

\begin{proof} First note that $\| [P_k, T] \| = \| [P_k, T^*] \| =
|\alpha_{i_k}| \to 0$ as $k \to \infty$.  Using the estimates $\| [P,
T^j] \| \leq \| PT^j - TPT^{j-1} \| + \| TPT^{j-1} - T^jP \| \leq \|
[P, T] \| \| T\|^{j-1} + \|T\|\|[P, T^{j-1} \|$, a simple induction
argument shows that $$\| [P_k, T^j] \| \leq j\|T\|^{j-1} \| [P_k, T]
\| = j\|T\|^{j-1}|\alpha_{i_k}|.$$ Since $S$ is a band operator
whose lower diagonals are determined by the powers of $T$ and upper
diagonals are determined by the powers of $T^*$, we may use the fact
that $\|[P, S]\| = \max \{ \|PSP^{\perp}\|, \| P^{\perp}SP\|\}$ to
arrive at the following inequality: $$\|[P_k, S]\| \leq M
q^{\prime}(\|T\|)|\alpha_{i_k}|.$$ With this inequality in hand, the
various statements above follow from the corresponding results in
sections 3 and 4.
\end{proof}

\begin{rem} Note that the rank of $P_k$ is $i_k$ and hence $P_kSP_k$
is an $i_k \times i_k$ matrix.  Thus the numerical efficiency in this
setting is directly tied to the rate of growth of $\{i_k\}$ and the
rate of decay of $\{ \alpha_{i_k} \}$ (as functions of $k$).
\end{rem}

\begin{rem} Evidently the theorem above can be extended with no
difficulty to the case of bilateral shifts on $l^2({\mathbb Z})$ where
one assumes the existence of subsequences tending to zero both in the
positive and negative directions.
\end{rem}

\begin{rem} Of course, the numerical ranges of $P_kSP_k$ also
converge to the numerical range of $S$ since the proof above shows
that $\|[P_k, S]\| \to 0$.  Note also that a similar result can be
worked out for any element in $C^*(T)$ and (complicated) error
estimates can be explicitly given for any polynomial in $T$ and $T^*$.
\end{rem}

\begin{rem} Let $f \in C({\mathbb T})$ be a continuous function on the
unit circle and $T_f$ be the corresponding Toeplitz operator on
$l^2({\mathbb N})$.  Of course, $T_f$ is rarely a band operator and
the diagonals of the matrix of $T_f$ never have subsequences tending
to zero (since they are constant).  However, it follows from BDF
theory that $T_f$ is quasidiagonal if and only if the winding number
of the curve $f({\mathbb T})$ is zero about every point in the
complement of $f({\mathbb T})$ (cf.\ \cite[Sections V.1 and
IX.7]{davidson}).  In this case there are asymptotically commuting
finite rank projections out there and we feel it would be interesting
(but perhaps hard?) to find them.  
\end{rem}

We conclude this section with an illustration of the usefulness of
this method.  Namely, we will give a new proof of the fact that the
spectra of quasidiagonal weighted shifts is as large as possible.  It
is well known that if $T$ is any weighted shift and $\rho(T)$ is the
spectral radius of $T$ then $\sigma(T) = \{ \lambda \in {\mathbb C}:
|\lambda| \leq \rho(T) \}$.  While our technology will only handle the
quasidiagonal case, the argument is simple and cute so we include a
proof.  It will be convenient to note that $$\rho(T) = \lim_{l \to
\infty} \|T^l\|^{1/l} = \lim_{l \to \infty} \bigg( \sup_{m \in
{\mathbb N}} |\alpha_m \alpha_{m+1}\cdots\alpha_{m+l}| \bigg)^{1/l},$$
whenever $T$ is a weighted shift with weight sequence $(\alpha_1,
\alpha_2, \ldots)$.

\begin{prop} 
Let $T$ be the weighted shift described at the beginning of this
section.  Then $\sigma(T) = \{ \lambda \in {\mathbb C}: |\lambda| \leq
\rho(T) \}$.
\end{prop}  

\begin{proof} Evidently we may assume that $\| T \| \leq 1$ by
rescaling, if necessary.

Since the spectrum of $T$ can be no larger than claimed, it suffices,
thanks to part (2) of Theorem \ref{thm:existspectrum}, to show that if
$$|\lambda| < \lim_{l \to \infty} \bigg( \sup_{m \in {\mathbb N}}
|\alpha_m \alpha_{m+1}\cdots\alpha_{m+l}| \bigg)^{1/l}$$ then for every
$\varepsilon > 0$, $\lambda$ belongs to the
$\varepsilon$-pseudospectra of $P_k T P_k$ for all sufficiently large
$k$. In fact, we will show a bit more (as it requires no extra effort
and simplifies notation).  Let $P_n$ be the projection onto the span
of the first $n$ basis vectors of $l^2({\mathbb N})$ and we will show
that $\lambda$ belongs to the $\varepsilon$-pseudospectra of $P_n T
P_n$ for all sufficiently large $n$.

First note that for all non-zero $\lambda$ the matrix $\lambda - P_n T
P_n$ is invertible.  A straightforward computation shows that its inverse is 
given by the following matrix. 
$$\begin{pmatrix}
\frac{1}{\lambda} & 0 & 0 & 0 & \cdots & 0 & 0\\
\frac{\alpha_1}{\lambda^2} & \frac{1}{\lambda} & 0 & 0 & \cdots 
& 0 & 0\\
\frac{\alpha_1\alpha_2}{\lambda^3} & \frac{\alpha_2}{\lambda^2} 
& \frac{1}{\lambda} & 0 & \cdots & 0 & 0\\
\frac{\alpha_1\alpha_2\alpha_3}{\lambda^4} & 
\frac{\alpha_2\alpha_3}{\lambda^3} & \frac{\alpha_3}{\lambda^2} 
& \frac{1}{\lambda} & \cdots & 0 & 0\\
\vdots & \vdots & \vdots & \vdots & \ddots & \vdots & \vdots\\ 
\frac{\prod_{1}^{n-1}\alpha_i}{\lambda^n} & 
\frac{\prod_{2}^{n-1}\alpha_i}{\lambda^{n-1}} & 
\frac{\prod_{3}^{n-1}\alpha_i}{\lambda^{n-2}} &
\frac{\prod_{4}^{n-1}\alpha_i}{\lambda^{n-3}} & \cdots & 
\frac{\alpha_{n -1}}{\lambda^2} & \frac{1}{\lambda} 
\end{pmatrix}$$

Our goal is to show that if $0 < |\lambda| < \lim_{l \to \infty}
\bigg( \sup_{m \in {\mathbb N}} |\alpha_m
\alpha_{m+1}\cdots\alpha_{m+l}| \bigg)^{1/l} \leq 1$ then for every
$\varepsilon > 0$ the norms of the matrices above are eventually
bigger than $1/ \varepsilon$.  To prove this it suffices to show that
the modulus of some entry in the matrices above is greater than $1/
\varepsilon$.  This, however, is the case because if $|\lambda| <
|\alpha_m \alpha_{m+1}\cdots\alpha_{m+l}|^{1/l} - \delta$ then $$(1 +
\delta)^l < \frac{|\alpha_m
\alpha_{m+1}\cdots\alpha_{m+l}|}{|\lambda|^l} < \frac{|\alpha_m
\alpha_{m+1}\cdots\alpha_{m+l}|}{|\lambda|^{l+1}}.$$ Hence, if
$\varepsilon > 0$ is given, we choose $l$ big enough that $(1 +
\delta)^l > 1/ \varepsilon$ and then we will be able to find a string
$\alpha_m, \ldots, \alpha_{m+l}$ such that $$\frac{1}{\varepsilon} <
\frac{|\alpha_m \alpha_{m+1}\cdots\alpha_{m+l}|}{|\lambda|^{l+1}}.$$
However, the numbers $$\frac{\alpha_m
\alpha_{m+1}\cdots\alpha_{m+l}}{\lambda^{l+1}}$$ are entries in all of
the matrices above (for sufficiently large $n$).
\end{proof}

\begin{rem} Note that the proof above adapts easily to bi-lateral
weighted shifts which have subsequences tending to zero in both the
positive and negative directions.
\end{rem}

\begin{rem}  
\label{thm:remark}
If $T$ denotes the bilateral shift on $l^2({\mathbb Z})$ and $P_n$
denotes the orthogonal projection onto the span of $\{ \delta_{-n},
\ldots, \delta_n \} \subset l^2({\mathbb Z})$ then the finite sections
$P_n T P_n$ are all nilpotent matrices and hence have spectrum $\{0\}$
for all $n$ (which, in particular, do not converge to the spectrum of
$T$).  Moreover, the proof above shows that for every $\varepsilon >
0$, $\liminf \sigma^{(\varepsilon)} (P_n TP_n)$ contains the entire
unit disc (see \cite[Definition 3.1]{HRS}).  Hence we can't recover
the spectrum of $T$ (i.e.\ the unit circle) by intersecting over
$\varepsilon > 0$ as in part (2) of Theorem \ref{thm:existspectrum}.
\end{rem}

%%%%%%%%%%%%%%%%%%%%%%%%%%%%%%%%%%%%%%%%%%%%%%%%%%%%%%%%%%%%%%%%%%
\section{Berg's Technique and Bilateral Band Operators}

In our final section we consider bilateral band operators and apply a
technique of Berg to produce explicit finite sections with the right
quasidiagonality requirements.  Berg's technique is a well-known, and
highly successful, way of taking a shift operator and approximating it
by cyclic shifts of a shorter length.  We will not attempt to explain
the philosophy behind this technique (as this is well done in
\cite[Section VI.4]{davidson}, for example) but rather will just write
down how it goes and check that it works.  The diligent reader who
actually checks all of the details will likely see what is going on in
the process.

In \cite{smucker} it was shown that a bilateral weighted shift is
quasidiagonal if and only if the weight sequence is ``block
balanced''.  The definition below is nothing but the obvious extension
of Smucker's notion to subsets of $l^{\infty}({\mathbb Z})$.

\begin{defn} A set ${\mathcal S} \subset l^{\infty}({\mathbb Z})$ is
called {\em block balanced} if for each $\varepsilon > 0$ and $k \in
{\mathbb N}$ there exist integers $m < 0$ and $n > 0$ such that $m + k
\leq 0$ and for each sequence $(a_i) \in {\mathcal S}$ we have $| a_{m+j}
- a_{n+j} | < \varepsilon$ for $0 \leq j \leq k$.
\end{defn}

This definition just says that we can find strings of arbitrary length
which are arbitrarily close to each other in the positive and negative
directions.  

We now describe the types of bilateral band operators to which the
quasidiagonal approach will apply.  Let $U$ be the bilateral shift on
$l^2 ({\mathbb Z})$ (i.e.\ $U(\delta_i) = \delta_{i+1}$ for the
canonical basis $\{ \delta_i \}_{i \in {\mathbb Z}}$ of $l^2 ({\mathbb
Z})$).  Assume $\{ f_{-p}, \ldots, f_{-1}, f_{0}, f_1, \ldots, f_p\}
\subset l^{\infty}({\mathbb Z})$ is a block balanced set and consider
the band operator $$T = \sum_{l = -p}^{p} f_l U^l.$$

Since we have assumed $\{ f_{-p}, \ldots, f_{-1}, f_{0}, f_1, \ldots,
f_p\}$ to be block balanced we can, for any sequence $\varepsilon_k
\to 0$ and natural numbers $s_k \to \infty$, find integers $m_k, n_k$
such that $m_k + s_k \leq 0 < n_k$ and $$| f_l (m_k + j) - f_l(n_k +
j) | < \varepsilon_k,$$ for all $-p \leq l \leq p$ and $0 \leq j \leq
s_k$.  Define constants $\alpha_{j,k} = \frac{1}{2}(1 -
\exp{\frac{(j+1)\pi}{s_k+2}})$ and $\beta_{j,k} = \frac{1}{2}(1 +
\exp{\frac{(j+1)\pi}{s_k+2}})$ for $0 \leq j \leq s_k$.  A calculation
shows that the set of vectors $\{\alpha_{0,k} \delta_{m_k} +
\beta_{0,k} \delta_{n_k}, \ldots, \alpha_{s_k,k} \delta_{m_k + s_k} +
\beta_{s_k,k} \delta_{n_k + s_k} \}$ is an orthonormal set.  Hence the
set
$$V_k = \{ \delta_i : m_k + s_k < i < n_k \} \cup \{\alpha_{0,k}
\delta_{m_k} + \beta_{0,k} \delta_{n_k}, \ldots, \alpha_{s_k,k}
\delta_{m_k + s_k} + \beta_{s_k,k} \delta_{n_k + s_k} \}$$ is also
orthonormal. 

Though it is not obvious, it turns out that the orthogonal projection
onto the subspace spanned by $V_k$ almost commutes with $T$. We will
denote this projection by $P_k$.

\begin{lem} 
\label{thm:lemma}
Keeping notation as above we have $$\| [T, P_k] \| \leq
\frac{(2p + 1)\varepsilon_k}{2} + \frac{\pi}{s_k + 2} \sum_{l = -p}^p
|l| \| f_l \|.$$
\end{lem} 

\begin{proof} Our first task will be to write down the 5x5 block
matrix of $P_k$ with respect to the decomposition $l^2 ({\mathbb Z}) =
H_1 \oplus H_2 \oplus H_3 \oplus H_4 \oplus H_5$, where $H_1 = span\{
\delta_j: j<m_k \}$, $H_2 = span\{ \delta_j: m_k \leq j \leq m_k +
s_k\}$, $H_3 = span\{ \delta_j: m_k + s_k < j < n_k\}$, $H_4 = span\{
\delta_j: n_k \leq j \leq n_k + s_k \}$ and $H_5 = span\{ \delta_j: n_k
+ s_k < j \}$.  

A tedious but straightforward calculation shows that in this
decomposition $P_k$ becomes
$$\begin{pmatrix}
0 & 0 & 0 & 0 & 0\\
0 & A_k & 0 & C_k & 0\\
0 & 0 & I & 0 & 0\\ 
0 & C_k^* & 0 & B_k & 0\\ 
0 & 0 & 0 & 0 & 0
\end{pmatrix}$$
where $I$ is the identity matrix, 
$$A_k = 
\begin{pmatrix}
|\alpha_0|^2 & 0 & \cdots & 0\\
0 & |\alpha_1|^2 & \cdots & 0\\
\vdots & \vdots & \ddots & \vdots\\ 
0 & 0 & \cdots & |\alpha_{s_k}|^2
\end{pmatrix}, 
B_k = 
\begin{pmatrix}
|\beta_0|^2 & 0 & \cdots & 0\\
0 & |\beta_1|^2 & \cdots & 0\\
\vdots & \vdots & \ddots & \vdots\\ 
0 & 0 & \cdots & |\beta_{s_k}|^2
\end{pmatrix}$$ and  
$$C_k = 
\begin{pmatrix}
\overline{\beta_0}\alpha_0 & 0 & \cdots & 0\\
0 & \overline{\beta_1}\alpha_1 & \cdots & 0\\
\vdots & \vdots & \ddots & \vdots\\ 
0 & 0 & \cdots & \overline{\beta_{s_k}}\alpha_{s_k}
\end{pmatrix},$$
where $\overline{\beta_j}$ is the complex conjugate of $\beta_j$.
Since the matrix of each $f_l$ is also block diagonal with respect to
the above decomposition we may write 
$$f_l = 
\begin{pmatrix}
F_1^l & 0 & 0 & 0 & 0\\
0 & F_2^l & 0 & 0 & 0\\
0 & 0 & F_3^l & 0 & 0\\ 
0 & 0 & 0 & F_4^l & 0\\ 
0 & 0 & 0 & 0 & F_5^l
\end{pmatrix}.$$
Hence $$\| [P_k, f_l] \| = \max\{\|F_2^l C_k - F_4^lC_k \|, \|F_2^l
C_k^* - F_4^lC_k^* \|\} \leq \varepsilon_k
\max\{|\overline{\beta_0}\alpha_0|, \ldots,
|\overline{\beta_{s_k}}\alpha_{s_k}|\} \leq \frac{\varepsilon_k}{2}.$$

Unfortunately, the matrix of the bilateral shift $U$ is a bit messy
with respect to the decomposition we have used above.  However,
conjugating by $U$ (i.e.\ considering $U \cdot U^*$) just has the
affect of shifting the matrix down the diagonal one entry.  This
remark together with the fact that $\| [U, P_k ] \| = \| U P_k U^* -
P_k \|$ allows one to estimate $\| [U, P_k ] \|$ above by $$\max\{
\big| |\alpha_j|^2 - |\alpha_{j+1}|^2\big|, \big| |\beta_j|^2 -
|\beta_{j+1}|^2\big|\} + \max\{\big| \overline{\beta_j}\alpha_j -
\overline{\beta_{j+1}}\alpha_{j+1}\big|\}.$$ But $|\alpha_j|^2 =
\frac{1}{2}(1- \cos(\frac{j\pi}{s_k+2}))$, $|\beta_j|^2 = \frac{1}{2}(1 +
\cos(\frac{j\pi}{s_k+2}))$ and $\overline{\beta_j}\alpha_j =
-\frac{i}{2} \sin(\frac{j\pi}{s_k+2})$.  Hence we get the inequality
$$\| [U, P_k ] \| \leq \frac{\pi}{s_k+2}.$$ 

Having the bounds on $\| [P_k, f_l] \|$ and $\| [U, P_k ] \|$ one
completes the proof by a standard interpolation argument, similar to
the proof of Theorem \ref{thm:unilateral}.
\end{proof}

\begin{rem} We leave it to the reader to check that $P_k (\delta_i)
\to \delta_i$, as $k \to \infty$, for each $i \in {\mathbb Z}$.  It
will not always be the case that $P_k \leq P_{k+1}$ but the fact that
$\|P_k(v) - v\| \to 0$ is good enough for all of the results in
sections 3, 4 and 5 to hold (i.e.\ we really don't have to have an
honest filtration in any of the proofs).
\end{rem}

\begin{rem} Note that the rank of $P_k$ is $n_k - m_k$ and hence the
numerical efficiency of this scheme will depend on how far out one
must go in order to find long strings which are close together.
\end{rem}

A number of interesting examples naturally give rise to tri-diagonal
operators on $l^2({\mathbb Z})$.  For example, a natural
discretization of $\frac{d}{dx}$ is $\frac{1}{2h}(U - U^*)$, where $h$
is a fixed step size.  A natural discretization of the one-dimensional
differential operator $(-\Delta + c)f = \frac{d^2}{dx^2}f + cf$, for
some constant $c$, would be $\frac{-1}{4h^2}(U + U^*) + (c +
\frac{1}{2h^2})I$ where $I \in l^{\infty}({\mathbb Z})$ is the
constant function $I(n) = 1$ for all $n$.  Also, Arveson's method of
discretizing Schr\"odinger operators of one-dimensional quantum
systems (cf.\ \cite{arveson}) naturally leads to an analysis of
operators of the form $U + U^* + f$ where $f(n) = v(\cos(n\theta))$
for some continuous function $v$ on the interval $[-1,1]$ and some
fixed positive number $\theta$.  Note that if $\frac{\theta}{\pi} \in
{\mathbb Q}$ then $f$ will be a {\em periodic} function while in
general it will be ``almost'' periodic and hence block balanced.
Hopefully these examples justify spending the remainder of this paper
analyzing how to implement Berg's technique on the class of
tri-diagonal operators with periodic main diagonal and constant
diagonals above and below the main diagonal.  The interested reader
should have no trouble now formulating a general result based on the
previous lemma and the results of section 4.  However, while this
method can be horribly inefficient in some cases, it works very well
in the periodic tri-diagonal case. Moreover, for those interested in
discretized Hamiltonians with $\frac{\theta}{\pi} \notin {\mathbb Q}$
we will describe how to use continued fractions to find good
``blocks'' at the very end.  (See also \cite{brown:PV} for improved
spectral computations in this case.)

For the remainder of this paper $T$ will be an operator of the form
$$T = c_1 U + c_{-1}U^* + f,$$ where $c_1, c_{-1} \in {\mathbb C}$ and
$f \in l^{\infty}({\mathbb Z})$ is a periodic function with period $p$
(i.e.\ $f(n + p) = f(n)$, for all $n \in {\mathbb Z}$).  Due to the
periodicity of $f$ and the constant coefficients on $U$ and $U^*$ we
may take $\varepsilon_k = 0$, $s_k = pk - 1$, $m_k = -s_k$ and $n_k =
1$ (notation as in the set-up for Lemma \ref{thm:lemma}).  Note that
with these choices the ``twisted'' basis used in Berg's technique is
$$V_k = \{\alpha_{0,k} \delta_{m_k} + \beta_{0,k} \delta_{1}, \ldots,
\alpha_{s_k,k} \delta_{0} + \beta_{s_k,k} \delta_{s_k} \}.$$ A tedious
but straightforward computation shows that if $P_k$ is the orthogonal
projection onto the span of $V_k$ then the matrix of $P_k T P_k$ with
respect to Berg's basis $V_k$ is
$$P_k T P_k = 
\begin{pmatrix}
\sigma_{0,0} & \sigma_{0,1} & 0 & 0 & \cdots & 0 & \sigma_{0,s_k}\\ 
\sigma_{1,0} & \sigma_{1,1} & \sigma_{1,2} & 0 & \cdots & 0 & 0\\ 
0 & \sigma_{2,1} & \sigma_{2,2} & \sigma_{2,3} & \cdots & 0 & 0\\
0 & 0 & \sigma_{3,2} & \sigma_{3,3} & \cdots  & 0 & 0\\ 
\vdots & \vdots & \vdots & \vdots & \ddots & \vdots & \vdots\\ 
0 & 0 & 0 & 0 & \cdots & \sigma_{s_k -1, s_k -1} & \sigma_{s_k -1, s_k}\\ 
\sigma_{s_k,0} & 0 & 0 & 0 & \cdots & \sigma_{s_k,s_k-1} & \sigma_{s_k,s_k}
\end{pmatrix},$$
where 
\begin{equation}
\begin{split}
&\sigma_{j,j} = |\alpha_{j,k}|^2 f(m_k+j) + |\beta_{j,k}|^2 f(1+j),\\
& \sigma_{j,j+1} = (\alpha_{j+1,k}\overline{\alpha_{j,k}} +
\beta_{j+1,k}\overline{\beta_{j,k}})c_{-1},\\
&\sigma_{j+1,j} =
(\alpha_{j,k}\overline{\alpha_{j+1,k}} +
\beta_{j,k}\overline{\beta_{j+1,k}})c_{1},\\
&\sigma_{0,s_k} =
\overline{\beta_{0,k}}\alpha_{s_k,k}c_1,\\
&\sigma_{s_k,0} =
\beta_{0,k}\overline{\alpha_{s_k,k}}c_{-1}.
\end{split}
\end{equation}

While it is unfortunate that the entries of the matrix $P_k TP_k$ are
more complicated than the entries of $T$ it is a significant feature
of Berg's technique that the property of being tri-diagonal (hence
sparse) is almost preserved.  Also note that $\| [T,P_k] \| \leq
\frac{2\pi}{pk+1}$ in this setting while the rank of $P_k$ is $pk$.
Hence the error estimates are inversely proportional to the rank of
the ``almost'' tri-diagonal matrix $P_kTP_k$.

We end this paper by putting all the error estimates together in the
case we are considering. The proof is immediate from Lemma
\ref{thm:lemma} and the results of sections 3, 4 and 5.

\begin{thm} Let $T, \ P_k \in B(l^2({\mathbb Z}))$ be as above. 
\begin{enumerate}
\item Let $y \in l^2({\mathbb Z})$ be given and assume that $T$ is
invertible.  If $k > \frac{1}{p}(2\pi\|T^{-1}\| - 1)$ then the
``almost'' tri-diagonal matrix $P_kTP_k$ is invertible and if $x_k$ is
the unique solution to $P_kTP_k x_k = P_k y$ then $$\| x - x_k \| \leq
\|T^{-1}\|(\| y - P_k y\| + \frac{2\pi \|x_k\|}{pk + 1}),$$ where $x$
is the unique solution to $Tx = y$.

\item For every $\varepsilon > 0$, $\sigma^{(\varepsilon)}(P_k TP_k) \to
\sigma^{(\varepsilon)}(T)$ and, moreover, $\sigma^{(\varepsilon)}(P_k
TP_k) \subset \sigma^{(\varepsilon + \frac{2\pi}{pk+1})}(T)$.

\item $\sigma_2(P_k TP_k) \to \sigma_2(T)$ and $\sigma_2(P_k TP_k)
\subset^{\sqrt{\frac{4\pi\|T\|}{pk+1}}} \sigma_2(T)$.

\item If $T = T^*$ then $\sigma(P_kTP_k) \to \sigma(T)$ and
$\sigma(P_kTP_k) \subset^{\frac{2\pi}{pk+1}} \sigma(T)$.

\end{enumerate}
\end{thm}

\begin{rem} The theory of continued fractions gives a remarkable
approximation property of irrational numbers (cf.\
\cite{hardy-wright}).  Namely, if $\alpha \in (0,1)$ is an irrational
number then there is a very simple and computable algorithm for
finding coprime integers $p_n, q_n$ such that $|\alpha -
\frac{p_n}{q_n}| \leq \frac{1}{q_n^2}$.  It follows that for any $0
\leq j \leq q_n$ we have the inequality $|j\alpha - j\frac{p_n}{q_n}|
\leq \frac{j}{q_n^2} \leq \frac{1}{q_n}$.  Hence we can apply Berg's
technique as above to a periodic approximation to a function of the
form $f(m) = \cos(m\theta)$ and still have good control on the error
estimates.
\end{rem}

\bibliographystyle{amsplain}

\begin{thebibliography}{10}

\bibitem{arveson} W.~ Arveson, \emph{The role of C$^*$-algebras in
infinite dimensional numerical linear algebra}, C$^*$-algebras:
1943-1993 (San Antonio, TX, 1993), 114--129, Contemp.\ Math.\
\textbf{167}, Amer.\ Math.\ Soc., Providence, RI, 1994.


%\bibitem{boca:book} F.P. Boca, \emph{Rotation C$^*$-algebras and
%almost Mathieu operators}, Theta Series in Advanced Mathematics, 1.
%The Theta Foundation, Bucharest, 2001.

\bibitem{bottcher} A.~ B\"{o}ttcher, \emph{C$^*$-algebras in numerical
analysis}, Irish Math.\ Soc.\ Bull.\ No.\ 45 (2000), 57--133.

\bibitem{brown:QDsurvey} N.P. Brown, \emph{On quasidiagonal
C$^*$-algebras}, math.OA/0008181, preprint. 

\bibitem{brown:AFD} N.P. Brown, \emph{Invariant means and finite
representation theory of C$^*$-algebras}, math.OA/0304009, preprint.

\bibitem{brown:PV} N.P. Brown, \emph{AF embeddings and the numerical
computation of spectra in irrational rotation algebras}, preprint.

\bibitem{brown:Herrero} N.P. Brown, \emph{Herrero's approximation
problem for quasidiagonal operators}, J.\ Funct.\ Anal.\ \textbf{186}
(2001), 360--365.

\bibitem{dadarlat} M.~ Dadarlat, \emph{On the approximation of
quasidiagonal C$^*$-algebras}, J.\ Funct.\ Anal.\ \textbf{167} (1999),
69--78.



\bibitem{davidson} K.R. Davidson, \emph{C$^*$-algebras by example},
Fields Institute Monographs 6, American Mathematical Society,
Providence, RI, 1996.

\bibitem{HRS} R.~ Hagen, S.~ Roch and B.~ Silbermann,
\emph{C$^*$-algebras and numerical analysis}. Monographs and Textbooks
in Pure and Applied Mathematics, 236. Marcel Dekker, Inc., New York,
2001.

\bibitem{hardy-wright} G.H. Hardy and E.M. Wright, \emph{An
introduction to the theory of numbers.} Fifth edition. The Clarendon
Press, Oxford University Press, New York, 1979.


\bibitem{halmos} P.R. Halmos, \emph{Ten problems in Hilbert space},
Bull.\ Amer.\ Math.\ Soc.\ \textbf{76} (1970), 887--933.

\bibitem{kirchberg} E.~ Kirchberg, \emph{Exact C$^*$-algebras, tensor
products and the classification of purely infinite algebras},
Proceedings of the International Congress of Mathematicians, Vol.\ 1,2
(Zurich, 1994), 943--954.

\bibitem{narayan} S.~Narayan, \emph{Quasidiagonality of direct sums of
weighted shifts}, Trans.\ Amer.\ Math.\ Soc.\ \textbf{332} (1992), 757
-- 774.

\bibitem{smucker} R.~Smucker, \emph{Quasidiagonal weighted shifts},
Pacific J.\ Math. \textbf{98} (1982), 173 -- 182.

\bibitem{szarek} S.~ Szarek, \emph{An exotic quasidiagonal operator},
J.\ Funct.\ Anal.\ \textbf{89} (1990), 274--290.

\bibitem{dvv:QDsurvey} D.V. Voiculescu, \emph{Around quasidiagonal
operators}, Integr. Equ. and Op. Thy. \textbf{17} (1993), 137--149.

\bibitem{voiculescu} D.V. Voiculescu, \emph{A note on quasidiagonal
operators}, Operator Theory: Advances and Applications, Vol.\ 32,
Birkhauser Verlag, Basel, 1988, 265--274.

\bibitem{wassermann} S.~ Wassermann, \emph{Exact C$^*$-algebras and
related topics}, Lecture Notes Series no.19, GARC, Seoul National
University, 1994. 


\end{thebibliography}

\providecommand{\bysame}{\leavevmode\hbox to3em{\hrulefill}\thinspace}

\end{document}